\documentclass[10pt,a4paper]{article}
\usepackage[left=20mm, right=15mm, top=15mm, bottom=15mm]{geometry}
\usepackage{subfig}
\usepackage{flushend}
\usepackage{indentfirst}
\usepackage{graphics}
\usepackage{amsmath}
\usepackage{graphicx}
\usepackage{epstopdf}
\usepackage{float}
\usepackage[font=footnotesize,labelfont=bf]{caption}
\usepackage[labelsep=period]{caption}

\def\be{\begin{equation}}
\def\ee{\end{equation}}
\def\bea{\begin{eqnarray}}
\def\eea{\end{eqnarray}}

\begin{document}

\title{{\huge \textbf{Analytical solutions of the Riccati equation with
coefficients satisfying integral or differential conditions with arbitrary
functions}}}
\date{}
\author{\textbf{Tiberiu Harko}$^{1,*}$, \textbf{Francisco S. N. Lobo}$^{2}$,
\textbf{M. K. Mak}$^{3}$ \\
\\
{\footnotesize $^{1}${Department of Mathematics, University College London,
Gower Street, London WC1E 6BT, United Kingdom}}\\
{\footnotesize $^{2}$Centro de Astronomia e Astrof\'{\i}sica da Universidade
de Lisboa, Campo Grande, Ed. C8 1749-016 Lisboa, Portugal}\\
{\footnotesize $^{3}$Department of Computing and Information Management,
Hong Kong Institute of Vocational Education,
Chai Wan,
Hong Kong, P. R. China%
}\\
{\footnotesize $^{*}$Corresponding Author: t.harko@ucl.ac.uk}}
\maketitle




\noindent \textbf{{\large {Abstract}}} \hspace{2pt} Ten new exact solutions
of the Riccati equation $dy/dx=a(x)+b(x)y+c(x)y^{2}$ are presented. The
solutions are obtained by assuming certain relations among the coefficients $%
a(x)$, $b(x)$ and $c(x)$ of the Riccati equation, in the form of some
integral or differential expressions, also involving some arbitrary
functions. By appropriately choosing the form of the coefficients of the
Riccati equation, with the help of the conditions imposed on the
coefficients, we obtain ten new integrability cases for the Riccati
equation. For each case the general solution of the Riccati equation is also
presented. The possibility of the application of the obtained mathematical
results for the study of anisotropic general relativistic stellar models is
also briefly considered.\newline

\newline

\noindent\hrulefill

\section{ Introduction}

The Riccati equation is one of the most studied first order non-linear
differential equations \cite{kamke, 1r, 2r}, and is given by
\begin{equation}
\frac{dy}{dx}=a(x)+b(x)y+c(x)y^{2},  \label{1}
\end{equation}%
where $a$, $b$, $c$ are arbitrary real functions of $x$, with $a,b,c\in
C^{\infty }(I)$, defined on a real interval $I\subseteq \Re $.

It is well-known that once a particular solution $y_{p}$ of the Riccati
equation is known, the general solution of Eq.~(\ref{1}) is given by
\begin{eqnarray}
y\left[ b\left( x\right) ,c\left( x\right) ,y_{p}\left( x\right) \right]
=y_{p}\left( x\right) + \frac{e^{\int^{x}\left[ b\left( \phi \right)
+2c\left( \phi \right) y_{p}\left( \phi \right) \right] d\phi }}{%
C-\int^{x}c\left( \psi \right) e^{\int^{\psi }\left[ b\left( \phi \right)
+2c\left( \phi \right) y_{p}\left( \phi \right) \right] d\phi }d\psi },
\label{GS}
\end{eqnarray}%
where $C$ is an arbitrary constant of integration and the particular
solution $y_{p}\left( x\right) $ satisfies the Riccati equation
\begin{equation}
\frac{dy_{p}}{dx}=a(x)+b(x)y_{p}\left( x\right) +c(x)y_{p}^{2}\left(
x\right) .  \label{GP}
\end{equation}%
If we know three particular solutions $y_{i}^{p}(x)$, $i=1,2,3$, then the
Riccati equation can be solved without quadratures \cite{kamke, 1r, 2r}.

Even when a particular solution is not known, the Riccati equation can be
integrated exactly if the coefficients of the equation satisfy some
conditions. For instance, if the coefficients of the Riccati equation
satisfy the following specific condition
\begin{equation}
a(x)+b(x)+c(x)\equiv 0,
\end{equation}%
then the solution of the Riccati equation is given by
\begin{equation}
y=\frac{K+\int {[c(x)+a(x)]E(x)dx}-E(x)}{K+\int {[c(x)+a(x)]E(x)dx}+E(x)},
\end{equation}
where $K$ is an arbitrary constant of integration, and $E(x)=\exp (\int {%
[c(x)-a(x)]dx}$ \cite{kamke}.

If the coefficients of the Riccati equation satisfy the more general
condition
\begin{equation}
\lambda ^{2}c(x)+\lambda \mu b(x)+\mu ^{2}a(x)\equiv 0,
\end{equation}%
where $\lambda $ and $\mu $ are arbitrary constants satisfying the condition
$|\lambda |+|\mu |>0$, then by means of the transformation $y(x)=\lambda
/\mu +u(x)$, the Riccati equation is transformed into a Bernoulli type
equation \cite{kamke},
\begin{equation}
u^{\prime }(x)=c(x)u^{2}(x)+\left[ \frac{2\lambda }{\mu }c(x)+b(x)\right]
u(x).
\end{equation}

Another interesting case is if the coefficients $a(x)$, $b(x)$, $c(x)$
satisfy the relation
\begin{eqnarray}
a(x)+b(x)+c(x)=\frac{d}{dx}\ln \left| \frac{\alpha (x)}{\beta (x)}\right| -
\frac{\alpha (x)-\beta (x)}{\alpha (x)\beta (x)}\left[ \alpha (x)c(x)-\beta
(x)a(x)\right] ,
\end{eqnarray}%
with $\alpha (x)$ and $\beta (x)$ properly chosen differentiable functions,
such that $\alpha \beta >0$, then the Riccati equation is integrable by
quadratures \cite{Strel}. If $c(x)\equiv 1$, and the functions $a(x)$ and $%
b(x)$ are polynomials satisfying the condition
\begin{equation}
\Delta =b^{2}(x)-2\frac{db(x)}{dx}-4a(x)\equiv \mathrm{constant},
\label{int1}
\end{equation}%
then
\begin{equation}
y_{\pm }(x)=-\frac{\left[ b(x)\pm \sqrt{\Delta }\right] }{2},  \label{int2}
\end{equation}%
%
%
%
%
%
%
are both solutions of the Riccati Eq.~(\ref{1}) \cite{kamke, 2r,Rain}. Very
recently, the integrability condition given by Eq.~(\ref{int1}) of the
Riccati Eq.~(\ref{1}), and of the reduced Riccati equation of the form $%
dy/dx=a(x)+c(x)y^{2}$, have been generalized in \cite{AMC,AMC1}.

Note that the Riccati equation plays a significant role in many fields of
applied and fundamental science \cite{R}. Some applications of the
integrability conditions for the case of the damped harmonic oscillator with
time dependent frequency and for solitonic wave have been briefly discussed
in \cite{AMC}. The applications of the integrability condition of the
reduced Riccati equation for the integration of the Schr\"{o}dinger and
Navier-Stokes equations have been also studied in \cite{AMC1}. The
integrability conditions derived in \cite{AMC,AMC1} have been applied in %
\cite{MH} to obtain a general solution of the Einstein's gravitational field
equations for the static spherically symmetric gravitational interior
space-time of an isotropic fluid sphere. The astrophysical analysis
indicates that this solution can be used as a realistic model for static
general relativistic high density objects, such as neutron stars. Riccati
equations also play an important role in cosmology. The mechanism of the
initial inflationary scenario of the Universe and of its late-time
acceleration can be described by assuming the existence of some
gravitationally coupled scalar fields $\phi $, with the inflaton field
generating inflation and the quintessence field being responsible for the
late accelerated expansion of the Universe. In the case of a scalar field
dominated Universe the dynamics of the expansion can be described by the
solutions of a Riccati type equation \cite{Scal}. In fact, the integrability
conditions of the Riccati equation, given by Eq. (\ref{1}), and its
applications to stellar and cosmological models have been extensively
discussed in the literature, and we refer the reader to \cite%
{MH6,MH3,MH4,MH5}.

Due to the nonlinear structure of the Riccati equation, the general solution
of the Riccati equation (\ref{1}), cannot be easily found. Therefore, one
has to use numerical techniques, or approximate method for obtaining its
solutions. Recently, a numerical method, using hybrid of block-pulse
functions, and Chebyshev polynomials for solving the Riccati equation has
been presented in \cite{CJB}. Various numerical methods such as Adomian
decomposition method, He's variational iteration method, homopoty
perturbation method, Taylor matrix method and Legendre wavelet method for
solving the Riccati equation have been proposed in \cite%
{TA,TS,AB,AB1,AB2,se,MM}.

It is the purpose of the present paper to present some further integrability
conditions of the Riccati equation, given by Eq.~(\ref{1}). The relations
among the coefficients of the Riccati equation involve some integral or
differential representations, as well as the presence of some arbitrary
functions. By appropriately choosing the form of the coefficients of the
Riccati equation we obtain ten new integrability conditions. For each case
the general solution of the Riccati equation is also obtained. The
possibility of the application of the obtained mathematical results for the
study of the anisotropic general relativistic stellar models is also briefly
considered.

The present paper is organized as follows. The ten new exact solutions of
the Riccati Eq.~(\ref{1}) are presented in Section~\ref{sect2}. Some
astrophysical applications of the solutions of the Riccati equation are
presented in Section \ref{sect3}. We conclude our results in Section~\ref%
{sect4}.

\section{Integrability conditions for the Riccati equation}
\label{sect2}

From the algebraic point of view the Riccati Eq.~(\ref{1}) is a
quadratic equation in $y\left( x\right) $. We consider that its
particular solutions $y_{\pm }^{p}(x)$ take the form
\begin{equation}
y_{\pm }^{p}(x)=\frac{-b\left( x\right) \pm \sqrt{%
b^{2}\left( x\right) -4a(x)c(x)+4c(x)\frac{dy^{p}}{dx}}}{2c(x)}.
\label{mn1}
\end{equation}

In order to obtain the general solution of the Riccati Eq.~(\ref{1}%
) with the help of Eq.~(\ref{mn1}), we introduce the new generating
function $f_{1}\left( x\right) $, satisfying the differential
condition given by
\be
b^{2}\left( x\right) +4c(x)\frac{dy^{p}}{dx}%
=f_{1}\left( x\right) ,
\ee
representing a first order differential equation in $y^p(x)$, and which can be immediately integrated to give the
particular solution of the Riccati Eq.~(\ref{1}). Therefore the general
solution of the Riccati Eq.~(\ref{1}) can be obtained through quadratures, since its particular
solution is known.

Based on the above algorithm, we shall show the detailed calculations in obtaining
Theorem 1 in Section 2.1, case 1. In order to make the paper readable and not
to have repetitive calculations, we shall not present the detailed
calculations for the rest of the theorems here. However, one may follow the
same procedures to obtain the Theorems 2-9 presented in this paper. In this
Section, we shall present  ten new integrability cases of
the Riccati Eq.~(\ref{1}). For each
case the general solution of the Riccati equation is explicitly obtained.

\subsection{Case 1: $a\left( x\right) =\frac{f_{1}\left( x\right) -\left\{
b\left( x\right) +c(x)\left[ \protect\int^{x}\frac{f_{1}\left( \protect\phi %
\right) -b^{2}\left( \protect\phi \right) }{2c(\protect\phi )}d\protect\phi %
-C_{1}\right] \right\} ^{2}}{4c(x)}$}

We assume that the arbitrary functions $b\left( x\right) $, $c(x)$ and $f_{1}(x)$ satisfy the differential
condition%
\begin{equation}
b^{2}\left( x\right) +4c(x)\frac{dy^{p}}{dx}=f_{1}\left( x\right) ,
\label{mn2}
\end{equation}
where we have introduced a generating function $f_{1}(x)\in
C^{\infty }(I)$ defined on a real interval $I\subseteq \Re $. The particular solutions $y_{\pm }^{p}\left( x\right) $
of the Riccati Eq.~(\ref{1}) take the form
\begin{equation}
y_{\pm }^{p}\left( x\right) \mathbf{=}\frac{-b\left( x\right) \pm
\sqrt{f_{1}\left( x\right) -4a(x)c(x)}}{2c(x)}\mathbf{=}\frac{1}{2}\left[
\int^{x}\frac{f_{1}\left( \phi \right) -b^{2}\left( \phi \right) }{2c(\phi )}%
d\phi -C_{1}\right] \mathbf{.}  \label{mn3}
\end{equation}
where we have used and integrated Eq.~(\ref{mn2}), and $C_{1}$%
 is an arbitrary constant of integration. By differentiating Eq.~(%
\ref{mn3}) with respect to $x$ yields the result
\begin{equation}
\frac{d}{dx}\left[ \frac{-b\left( x\right) \pm \sqrt{f_{1}\left( x\right)
-4a(x)c(x)}}{c(x)}\right] \mathbf{=}\frac{f_{1}\left( x\right) -b^{2}\left(
x\right) }{2c(x)}\mathbf{.}  \label{mn4}
\end{equation}

Eq.~(\ref{mn4}) can be integrated to give the coefficient $a\left(
x\right) $ of the Riccati Eq.~(\ref{1}) as
\begin{equation}
a\left( x\right) \mathbf{=}\frac{f_{1}\left( x\right) -\left\{
b\left( x\right) +c(x)\left[ \int^{x}\frac{f_{1}\left( \phi \right)
-b^{2}\left( \phi \right) }{2c(\phi )}d\phi -C_{1}\right] \right\} ^{2}}{%
4c(x)}\mathbf{,}  \label{a1}
\end{equation}
where  $C_{1}$is an arbitrary constant of integration. By substituting
Eq.~(\ref{a1}) into the Riccati Eq.~(\ref{1}), the latter can be expressed
as
\begin{equation}\label{TH1}
\frac{dy}{dx}=\frac{f_{1}\left( x\right) -\left\{ b\left( x\right) +c(x)%
\left[ \int^{x}\frac{f_{1}\left( \phi \right) -b^{2}\left( \phi \right) }{%
2c(\phi )}d\phi -C_{1}\right] \right\} ^{2}}{4c(x)}+b(x)y+c(x)y^{2}.
\end{equation}

Therefore we obtain the following:

\textbf{Theorem 1.} If the coefficient $a(x)$ of the Riccati Eq.~(%
\ref{1}) satisfies the integral condition (\ref{a1}), then the general
solution of the Riccati Eq.~(\ref{TH1}) is given by %
\begin{equation}
y(x)=\frac{e^{\int^{x}\left\{ b\left( \phi \right) +c(\phi )\left[
\int^{\phi }\frac{f_{1}\left( \psi \right) -b^{2}\left( \psi \right) }{%
2c(\psi )}d\psi -C_{1}\right] \right\} d\phi }}{C_{0}-\int^{x}c\left( \psi
\right) e^{\int^{\psi }\left\{ b\left( \phi \right) +c(\phi )\left[
\int^{\phi }\frac{f_{1}\left( \psi \right) -b^{2}\left( \psi \right) }{%
2c(\psi )}d\psi -C_{1}\right] \right\} d\phi }d\psi }+\frac{1}{2}\left[
\int^{x}\frac{f_{1}\left( \phi \right) -b^{2}\left( \phi \right) }{2c(\phi )}%
d\phi -C_{1}\right] ,
\end{equation}%
where $C_{0}$ is an arbitrary constant of integration.

\subsection{Case 2: $a\left( x\right) =\frac{d}{dx}\left[ \frac{-b\left(
x\right) \pm \protect\sqrt{f_{2}\left( x\right) +b^{2}(x)}}{2c(x)}\right] -%
\frac{f_{2}\left( x\right) }{4c(x)}$}

Next, we assume that the arbitrary functions $a\left( x\right) $ satisfies
the differential condition%
\begin{equation}
a\left( x\right) =\frac{d}{dx}\left[ \frac{-b\left( x\right) \pm \sqrt{%
f_{2}\left( x\right) +b^{2}(x)}}{2c(x)}\right] -\frac{f_{2}\left( x\right) }{%
4c(x)},  \label{d1}
\end{equation}
where we have introduced a new arbitrary function $f_{2}(x)\in C^{\infty
}(I) $ defined on a real interval $I\subseteq \Re $. By substituting Eq.~(%
\ref{d1}) into the Riccati Eq.~(\ref{1}), the latter can be expressed as
\begin{eqnarray}
\frac{dy_{\pm }}{dx}=\left\{ \frac{d}{dx}\left[ \frac{-b\left( x\right) \pm
\sqrt{f_{2}\left( x\right) +b^{2}(x)}}{2c(x)}\right] -\frac{f_{2}\left(
x\right) }{4c(x)}\right\} + b(x)y_{\pm }+c(x)y_{\pm }^{2}.  \label{TH2}
\end{eqnarray}

Therefore we obtain the following:

\textbf{Theorem 2}. If the coefficient $a(x)$ of the Riccati Eq.~(\ref{1})
satisfies the differential condition (\ref{d1}), then the general solutions
of the Riccati Eq.~(\ref{TH2}) are given by
\begin{eqnarray}
y_{\pm }(x)=\frac{e^{\pm \int^{x}\sqrt{f_{2}\left( \phi \right) +b^{2}\left(
\phi \right) }d\phi }}{C_{\pm 2}-\int^{x}c\left( \psi \right) e^{\pm
\int^{\psi }\sqrt{f_{2}\left( \phi \right) +b^{2}\left( \phi \right) }d\phi
}d\psi }+ \left[ \frac{-b\left( x\right) \pm \sqrt{f_{2}\left( x\right)
+b^{2}(x)}}{2c(x)}\right] ,
\end{eqnarray}
where $C_{\pm 2}$ are arbitrary constants of integration.

\subsection{\protect\bigskip Case 3: $b\left( x\right) =\frac{f_{2}\left(
x\right) -4c^{2}\left( x\right) \left\{ \protect\int^{x}\left[ a\left(
\protect\phi \right) +\frac{f_{2}\left( \protect\phi \right) }{4c\left(
\protect\phi \right) }\right] d\protect\phi -C_{3}\right\} ^{2}}{4c\left(
x\right) \left\{ \protect\int^{x}\left[ a\left( \protect\phi \right) +\frac{%
f_{2}\left( \protect\phi \right) }{4c\left( \protect\phi \right) }\right] d%
\protect\phi -C_{3}\right\} }$}

We now assume that the arbitrary function $b\left( x\right) $ satisfies the
integral condition%
\begin{equation}
b\left( x\right) =\frac{f_{2}\left( x\right) -4c^{2}\left( x\right) \left\{
\int^{x}\left[ a\left( \phi \right) +\frac{f_{2}\left( \phi \right) }{%
4c\left( \phi \right) }\right] d\phi -C_{3}\right\} ^{2}}{4c\left( x\right)
\left\{ \int^{x}\left[ a\left( \phi \right) +\frac{f_{2}\left( \phi \right)
}{4c\left( \phi \right) }\right] d\phi -C_{3}\right\} },  \label{b5a}
\end{equation}
where $C_{3}$ is an arbitrary constant. By substituting Eq.~(\ref{b5a}) into
the Riccati Eq.~(\ref{1}), the latter can be expressed as
\begin{eqnarray}
\frac{dy}{dx}=a\left( x\right) + \frac{f_{2}\left( x\right) -4c^{2}\left(
x\right) \left\{ \int^{x}\left[ a\left( \phi \right) +\frac{f_{2}\left( \phi
\right) }{4c\left( \phi \right) }\right] d\phi -C_{3}\right\} ^{2}}{4c\left(
x\right) \left\{ \int^{x}\left[ a\left( \phi \right) +\frac{f_{2}\left( \phi
\right) }{4c\left( \phi \right) }\right] d\phi -C_{3}\right\} }y+ c(x)y^{2}.
\label{TH4}
\end{eqnarray}

Therefore we obtain the following:

\textbf{Theorem 3}. If the coefficient $b(x)$ of the Riccati Eq.~(\ref{1})
satisfies the integral condition (\ref{b5a}), then the general solution of
the Riccati Eq.~(\ref{TH4}) is given by

\begin{eqnarray}
y\left( x\right) &=&\frac{e^{-\int^{x}\sqrt{f_{2}\left( \omega \right)
+\left\{ \frac{f_{2}\left( \omega \right) -4c^{2}\left( \omega \right)
\left\{ \int^{\omega }\left[ a\left( \phi \right) +\frac{f_{2}\left( \phi
\right) }{4c\left( \phi \right) }\right] d\phi -C_{3}\right\} ^{2}}{4c\left(
\omega \right) \left\{ \int^{\omega }\left[ a\left( \phi \right) +\frac{%
f_{2}\left( \phi \right) }{4c\left( \phi \right) }\right] d\phi
-C_{3}\right\} }\right\} ^{2}}d\omega }}{C_{4}-\int^{x}c\left( \psi \right)
e^{-\int^{\psi }\sqrt{f_{2}\left( \omega \right) +\left\{ \frac{f_{2}\left(
\omega \right) -4c^{2}\left( \omega \right) \left\{ \int^{\omega }\left[
a\left( \phi \right) +\frac{f_{2}\left( \phi \right) }{4c\left( \phi \right)
}\right] d\phi -C_{3}\right\} ^{2}}{4c\left( \omega \right) \left\{
\int^{\omega }\left[ a\left( \phi \right) +\frac{f_{2}\left( \phi \right) }{%
4c\left( \phi \right) }\right] d\phi -C_{3}\right\} }\right\} ^{2}}d\omega
}d\psi }  \notag \\
&&-\frac{1}{2c\left( x\right) }\left\{
\begin{array}{c}
\sqrt{f_{2}\left( x\right) +\left\{ \frac{f_{2}\left( x\right) -4c^{2}\left(
x\right) \left\{ \int^{x}\left[ a\left( \phi \right) +\frac{f_{2}\left( \phi
\right) }{4c\left( \phi \right) }\right] d\phi -C_{3}\right\} ^{2}}{4c\left(
x\right) \left\{ \int^{x}\left[ a\left( \phi \right) +\frac{f_{2}\left( \phi
\right) }{4c\left( \phi \right) }\right] d\phi -C_{3}\right\} }\right\} ^{2}}%
+ \\
\frac{f_{2}\left( x\right) -4c^{2}\left( x\right) \left\{ \int^{x}\left[
a\left( \phi \right) +\frac{f_{2}\left( \phi \right) }{4c\left( \phi \right)
}\right] d\phi -C_{3}\right\} ^{2}}{4c\left( x\right) \left\{ \int^{x}\left[
a\left( \phi \right) +\frac{f_{2}\left( \phi \right) }{4c\left( \phi \right)
}\right] d\phi -C_{3}\right\} }%
\end{array}%
\right\} ,
\end{eqnarray}%
where $C_{4}$ is an arbitrary constant of integration.

\subsection{ Case 4: $c\left( x\right) =\frac{\left[ -b\left( x\right) \pm
\protect\sqrt{f_{2}\left( x\right) +b^{2}(x)}\right] e^{-\frac{1}{2}\protect%
\int^{x}\frac{f_{2}\left( \protect\phi \right) }{-b\left( \protect\phi %
\right) \pm \protect\sqrt{f_{2}\left( \protect\phi \right) +b^{2}(\protect%
\phi )}}d\protect\phi }}{2\left[ C_{\pm 5}+\protect\int^{x}a\left( \protect%
\psi \right) e^{-\frac{1}{2}\protect\int^{\protect\psi }\frac{f_{2}\left(
\protect\phi \right) }{-b\left( \protect\phi \right) \pm \protect\sqrt{%
f_{2}\left( \protect\phi \right) +b^{2}(\protect\phi )}}d\protect\phi }d%
\protect\psi \right] }$}

We assume that the arbitrary function $c\left( x\right) $ satisfies the
integral condition%
\begin{equation}
c\left( x\right) =\frac{\left[ -b\left( x\right) \pm \sqrt{f_{2}\left(
x\right) +b^{2}(x)}\right] e^{-\frac{1}{2}\int^{x}\frac{f_{2}\left( \phi
\right) }{-b\left( \phi \right) \pm \sqrt{f_{2}\left( \phi \right)
+b^{2}(\phi )}}d\phi }}{2\left[ C_{\pm 5}+\int^{x}a\left( \psi \right) e^{-%
\frac{1}{2}\int^{\psi }\frac{f_{2}\left( \phi \right) }{-b\left( \phi
\right) \pm \sqrt{f_{2}\left( \phi \right) +b^{2}(\phi )}}d\phi }d\psi %
\right] },  \label{c1}
\end{equation}
where $C_{\pm 5}$ are arbitrary constants. By substituting Eq.~(\ref{c1})
into the Riccati Eq.~(\ref{1}), the latter can be expressed as
\begin{equation}
\frac{dy_{\pm }}{dx}=a\left( x\right) +b(x)y_{\pm }+\frac{\left[ -b\left(
x\right) \pm \sqrt{f_{2}\left( x\right) +b^{2}(x)}\right] e^{-\frac{1}{2}%
\int^{x}\frac{f_{2}\left( \phi \right) }{-b\left( \phi \right) \pm \sqrt{%
f_{2}\left( \phi \right) +b^{2}(\phi )}}d\phi }}{2\left[ C_{\pm
5}+\int^{x}a\left( \psi \right) e^{-\frac{1}{2}\int^{\psi }\frac{f_{2}\left(
\phi \right) }{-b\left( \phi \right) \pm \sqrt{f_{2}\left( \phi \right)
+b^{2}(\phi )}}d\phi }d\psi \right] }y_{\pm }^{2}.  \label{TH3}
\end{equation}

Therefore we obtain the following:

\textbf{Theorem 4}. If the coefficient $c(x)$ of the Riccati Eq.~(\ref{1})
satisfies the integral condition (\ref{c1}), then the general solutions of
the Riccati Eq.~(\ref{TH3}) are given by

\begin{eqnarray}
y_{\pm }\left( x\right) &=&\frac{e^{\pm \int^{x}\sqrt{f_{2}\left( \phi
\right) +b^{2}(\phi )}d\phi }}{C_{\pm 6}-\int^{x}\frac{\left[ -b\left(
\omega \right) \pm \sqrt{f_{2}\left( \omega \right) +b^{2}(\omega )}\right]
e^{-\frac{1}{2}\int^{\omega }\frac{f_{2}\left( \phi \right) }{-b\left( \phi
\right) \pm \sqrt{f_{2}\left( \phi \right) +b^{2}(\phi )}}d\phi }}{2\left[
C_{\pm 5}+\int^{\omega }a\left( \psi \right) e^{-\frac{1}{2}\int^{\psi }%
\frac{f_{2}\left( \phi \right) }{-b\left( \phi \right) \pm \sqrt{f_{2}\left(
\phi \right) +b^{2}(\phi )}}d\phi }d\psi \right] }e^{\pm \int^{\omega }\sqrt{%
f_{2}\left( \phi \right) +b^{2}(\phi )}d\phi }d\omega }  \notag \\
&&+\left[ C_{\pm 5}+\int^{x}a\left( \psi \right) e^{-\frac{1}{2}\int^{\psi }%
\frac{f_{2}\left( \phi \right) }{-b\left( \phi \right) \pm \sqrt{f_{2}\left(
\phi \right) +b^{2}(\phi )}}d\phi }d\psi \right] e^{\frac{1}{2}\int^{x}\frac{%
f_{2}\left( \phi \right) }{-b\left( \phi \right) \pm \sqrt{f_{2}\left( \phi
\right) +b^{2}(\phi )}}d\phi },
\end{eqnarray}

where $C_{\pm 6}$ are arbitrary constants of integration.

\subsection{Case 5: $a\left( x\right) =\frac{1}{4}\left\{ \frac{f_{3}\left(
x\right) }{c\left( x\right) }-2b\left( x\right) \left[ \protect\int^{x}\frac{%
f_{3}\left( \protect\phi \right) }{2c(\protect\phi )}d\protect\phi -C_{7}%
\right] -c\left( x\right) \left[ \protect\int^{x}\frac{f_{3}\left( \protect%
\phi \right) }{2c(\protect\phi )}d\protect\phi -C_{7}\right] ^{2}\right\} $}

Assume now that the arbitrary function $a(x)$ satisfies the integral
condition%
\begin{equation}
a\left( x\right) =\frac{1}{4}\left\{ \frac{f_{3}\left( x\right) }{c\left(
x\right) }-2b\left( x\right) \left[ \int^{x}\frac{f_{3}\left( \phi \right) }{%
2c(\phi )}d\phi -C_{7}\right] -c\left( x\right) \left[ \int^{x}\frac{%
f_{3}\left( \phi \right) }{2c(\phi )}d\phi -C_{7}\right] ^{2}\right\} ,
\label{d4}
\end{equation}
where $\ f_{3}(x)\in C^{\infty }(I)$ is an arbitrary function defined on a
real interval $I\subseteq \Re $ and $C_{7}$ is an arbitrary constant. By
substituting Eq.~(\ref{d4}) into the Riccati Eq.~(\ref{1}), the latter can
be expressed as

\begin{equation}
\frac{dy}{dx}=\frac{1}{4}\left\{ \frac{f_{3}\left( x\right) }{c\left(
x\right) }-2b\left( x\right) \left[ \int^{x}\frac{f_{3}\left( \phi \right) }{%
2c(\phi )}d\phi -C_{7}\right] -c\left( x\right) \left[ \int^{x}\frac{%
f_{3}\left( \phi \right) }{2c(\phi )}d\phi -C_{7}\right] ^{2}\right\}
+b(x)y+c(x)y^{2}.  \label{TH6}
\end{equation}

Therefore we obtain the following:

\textbf{Theorem 5}. If the coefficient $a(x)$ of the Riccati Eq.~(\ref{1})
satisfies the integral condition (\ref{d4}), then the general solution of
the Riccati Eq.~(\ref{TH6}) is given by

\begin{eqnarray}
y(x) &=&\frac{e^{\int^{x}\sqrt{b^{2}\left( \psi \right) +c\left( \psi
\right) \left[ \int^{\psi }\frac{f_{3}\left( \phi \right) }{2c(\phi )}d\phi
-C_{7}\right] \left\{ 2b\left( \psi \right) +c\left( \psi \right) \left[
\int^{\psi }\frac{f_{3}\left( \phi \right) }{2c(\phi )}d\phi -C_{7}\right]
\right\} }d\psi }}{C_{8}-\int^{x}c\left( \omega \right) e^{\int^{\omega }%
\sqrt{b^{2}\left( \psi \right) +c\left( \psi \right) \left[ \int^{\psi }%
\frac{f_{3}\left( \phi \right) }{2c(\phi )}d\phi -C_{7}\right] \left\{
2b\left( \psi \right) +c\left( \psi \right) \left[ \int^{\psi }\frac{%
f_{3}\left( \phi \right) }{2c(\phi )}d\phi -C_{7}\right] \right\} }d\psi
}d\omega }  \notag \\
&&+\frac{1}{2c\left( x\right) }\left\{ -b\left( x\right) +\sqrt{b^{2}\left(
x\right) +c\left( x\right) \left[ \int^{x}\frac{f_{3}\left( \phi \right) }{%
2c(\phi )}d\phi -C_{7}\right] \left\{ 2b\left( x\right) +c\left( x\right) %
\left[ \int^{x}\frac{f_{3}\left( \phi \right) }{2c(\phi )}d\phi -C_{7}\right]
\right\} }\right\} ,
\end{eqnarray}
where $C_{8}$ is an arbitrary constant of integration.

\subsection{ Case 6: $b\left( x\right) =\frac{f_{3}\left( x\right) -4a\left(
x\right) c\left( x\right) -c^{2}\left( x\right) \left[ \protect\int^{x}\frac{%
f_{3}\left( \protect\phi \right) }{2c(\protect\phi )}d\protect\phi -C_{7}%
\right] ^{2}}{2c\left( x\right) \left[ \protect\int^{x}\frac{f_{3}\left(
\protect\phi \right) }{2c(\protect\phi )}d\protect\phi -C_{7}\right] }$}

In this case, we assume that the coefficient $b\left( x\right) $ of the
Riccati Eq.~(\ref{1}) satisfies the integral condition
\begin{equation}
b\left( x\right) =\frac{f_{3}\left( x\right) -4a\left( x\right) c\left(
x\right) -c^{2}\left( x\right) \left[ \int^{x}\frac{f_{3}\left( \phi \right)
}{2c(\phi )}d\phi -C_{7}\right] ^{2}}{2c\left( x\right) \left[ \int^{x}\frac{%
f_{3}\left( \phi \right) }{2c(\phi )}d\phi -C_{7}\right] }.  \label{b7}
\end{equation}

By substituting Eq.~(\ref{b7}) into the Riccati Eq.~(\ref{1}), the latter
can be expressed as
\begin{equation}
\frac{dy}{dx}=a\left( x\right) +\frac{f_{3}\left( x\right) -4a\left(
x\right) c\left( x\right) -c^{2}\left( x\right) \left[ \int^{x}\frac{%
f_{3}\left( \phi \right) }{2c(\phi )}d\phi -C_{7}\right] ^{2}}{2c\left(
x\right) \left[ \int^{x}\frac{f_{3}\left( \phi \right) }{2c(\phi )}d\phi
-C_{7}\right] }y+c(x)y^{2}.  \label{TH7}
\end{equation}

Therefore we obtain the following:

\textbf{Theorem 6}. If the coefficient $b(x)$ of the Riccati Eq.~(\ref{1})
satisfies the integral condition (\ref{b7}), then the general solution of
the Riccati Eq.~(\ref{TH7}) is given by

\begin{eqnarray}
y(x) &=&\frac{e^{-\int^{x}\sqrt{\left\{ \frac{f_{3}\left( \psi \right)
-4a\left( \psi \right) c\left( \psi \right) -c^{2}\left( \psi \right) \left[
\int^{\psi }\frac{f_{3}\left( \phi \right) }{2c(\phi )}d\phi -C_{7}\right]
^{2}}{2c\left( \psi \right) \left[ \int^{\psi }\frac{f_{3}\left( \phi
\right) }{2c(\phi )}d\phi -C_{7}\right] }\right\} ^{2}-4a\left( \psi \right)
c\left( \psi \right) +f_{3}\left( \psi \right) }d\psi }}{C_{9}-\int^{x}c%
\left( \omega \right) e^{-\int^{\omega }\sqrt{\left\{ \frac{f_{3}\left( \psi
\right) -4a\left( \psi \right) c\left( \psi \right) -c^{2}\left( \psi
\right) \left[ \int^{\psi }\frac{f_{3}\left( \phi \right) }{2c(\phi )}d\phi
-C_{7}\right] ^{2}}{2c\left( \psi \right) \left[ \int^{\psi }\frac{%
f_{3}\left( \phi \right) }{2c(\phi )}d\phi -C_{7}\right] }\right\}
^{2}-4a\left( \psi \right) c\left( \psi \right) +f_{3}\left( \psi \right) }%
d\psi }d\omega }  \notag \\
&&+\frac{1}{2c\left( x\right) }\left\{
\begin{array}{c}
-\sqrt{\left\{ \frac{f_{3}\left( x\right) -4a\left( x\right) c\left(
x\right) -c^{2}\left( x\right) \left[ \int^{x}\frac{f_{3}\left( \phi \right)
}{2c(\phi )}d\phi -C_{7}\right] ^{2}}{2c\left( x\right) \left[ \int^{x}\frac{%
f_{3}\left( \phi \right) }{2c(\phi )}d\phi -C_{7}\right] }\right\}
^{2}-4a\left( x\right) c\left( x\right) +f_{3}\left( x\right) }+ \\
\frac{4a\left( x\right) c\left( x\right) -f_{3}\left( x\right) +c^{2}\left(
x\right) \left[ \int^{x}\frac{f_{3}\left( \phi \right) }{2c(\phi )}d\phi
-C_{7}\right] ^{2}}{2c\left( x\right) \left[ \int^{x}\frac{f_{3}\left( \phi
\right) }{2c(\phi )}d\phi -C_{7}\right] }%
\end{array}%
\right\} ,
\end{eqnarray}%
where $C_{9}$ is arbitrary constant of integration.

\subsection{Case 7: $\ a\left( x\right) =\frac{1}{4c\left( x\right) }\left[
2c\left( x\right) \frac{d}{dx}\left( \frac{f_{4}}{c}\right) -f_{
4}^{2}\left( x\right) -2b\left( x\right) f_{4}\left( x\right) \right] $}

We assume that the coefficient $a\left( x\right) $ of the Riccati Eq.~(\ref%
{1}) satisfies the differential condition%
\begin{equation}
a\left( x\right) =\frac{1}{4c\left( x\right) }\left[ 2c\left( x\right) \frac{%
d}{dx}\left( \frac{f_{4}}{c}\right) -f_{4}^{2}\left( x\right) -2b\left(
x\right) f_{4}\left( x\right) \right] ,  \label{d5}
\end{equation}%
where $f_{4}\left( x\right) \in C^{\infty }(I)$ is an arbitrary function
defined on a real interval $I\subseteq \Re $.

By substituting Eq.~(\ref{d5}) into the Riccati Eq.~(\ref{1}), the latter
can be expressed as
\begin{equation}
\frac{dy}{dx}=\frac{1}{4c\left( x\right) }\left[ 2c\left( x\right) \frac{d}{%
dx}\left( \frac{f_{4}}{c}\right) -f_{4}^{2}\left( x\right) -2b\left(
x\right) f_{4}\left( x\right) \right] +b(x)y+c(x)y^{2}.  \label{TH8}
\end{equation}

Therefore we obtain the following:

\textbf{Theorem 7}. If the coefficient $a(x)$ of the Riccati Eq.~(\ref{1})
satisfies the differential condition~(\ref{d5}), then the general solution
of the Riccati Eq.~(\ref{TH8}) is given by
\begin{equation}
y(x)=\frac{e^{\int^{x}\left[ b\left( \phi \right) +f_{4}\left( \phi \right) %
\right] d\phi }}{C_{10}-\int^{x}c\left( \psi \right) e^{\int^{\psi }\left[
b\left( \phi \right) +f_{4}\left( \phi \right) \right] d\phi }d\psi }+\frac{%
f_{4}\left( x\right) }{2c\left( x\right) },
\end{equation}%
where $C_{10}$ is an arbitrary constant of integration.

\subsection{ Case 8: $b\left( x\right) =\frac{1}{f_{4}\left( x\right) }\left[
c\left( x\right) \frac{d}{dx}\left( \frac{f_{4}}{c}\right) -\frac{%
f_{4}^{2}\left( x\right) }{2}-2a\left( x\right) c\left( x\right) \right] $}

We assume that the coefficient $b\left( x\right) $ of the Riccati Eq.~(\ref%
{1}) satisfies the differential condition
\begin{equation}
b\left( x\right) =\frac{1}{f_{4}\left( x\right) }\left[ c\left( x\right)
\frac{d}{dx}\left( \frac{f_{4}}{c}\right) -\frac{f_{4}^{2}\left( x\right) }{2%
}-2a\left( x\right) c\left( x\right) \right] .  \label{b8}
\end{equation}%
By substituting Eq.~(\ref{b8}) into the Riccati Eq.~(\ref{1}), the latter
can be expressed as
\begin{equation}
\frac{dy}{dx}=a\left( x\right) +\frac{1}{f_{4}\left( x\right) }\left[
c\left( x\right) \frac{d}{dx}\left( \frac{f_{4}}{c}\right) -\frac{%
f_{4}^{2}\left( x\right) }{2}-2a\left( x\right) c\left( x\right) \right]
y+c(x)y^{2}.  \label{TH8a}
\end{equation}

Therefore we obtain the following:

\textbf{Theorem 8}. If the coefficient $b(x)$ of the Riccati Eq.~(\ref{1})
satisfies the differential condition~(\ref{b8}), then the general solution
of the Riccati Eq.~(\ref{TH8a}) is given by
\begin{equation}
y(x)=\frac{e^{\int^{x}\left\{ \frac{1}{f_{4}\left( \phi \right) }\left[
c\left( \phi \right) \frac{d}{d\phi }\left( \frac{f_{4}}{c}\right) -\frac{%
f_{4}^{2}\left( \phi \right) }{2}-2a\left( \phi \right) c\left( \phi \right) %
\right] +f_{4}\left( \phi \right) \right\} d\phi }}{C_{11}-\int^{x}c\left(
\psi \right) e^{\int^{\psi }\left\{ \frac{1}{f_{4}\left( \phi \right) }\left[
c\left( \phi \right) \frac{d}{d\phi }\left( \frac{f_{4}}{c}\right) -\frac{%
f_{4}^{2}\left( \phi \right) }{2}-2a\left( \phi \right) c\left( \phi \right) %
\right] +f_{4}\left( \phi \right) \right\} d\phi }d\psi }+\frac{f_{4}\left(
x\right) }{2c\left( x\right) },
\end{equation}%
where $C_{11}$ is an arbitrary constant of integration.

\subsection{Case 9: $c\left( x\right) =\frac{f_{4}\left( x\right) e^{-%
\protect\int^{x}\left[ \frac{f_{4}\left( \protect\phi \right) }{2}+b\left(
\protect\phi \right) \right] d\protect\phi }}{C_{12}+2\protect\int%
^{x}a\left( \protect\psi \right) e^{-\protect\int^{\protect\psi }\left[
\frac{f_{4}\left( \protect\phi \right) }{2}+b\left( \protect\phi \right) %
\right] d\protect\phi }d\protect\psi }$}

We assume that the coefficient $c\left( x\right) $ of the Riccati Eq.~(\ref%
{1}) satisfies the integral condition
\begin{equation}
c\left( x\right) =\frac{f_{4}\left( x\right) e^{-\int^{x}\left[ \frac{%
f_{4}\left( \phi \right) }{2}+b\left( \phi \right) \right] d\phi }}{%
C_{12}+2\int^{x}a\left( \psi \right) e^{-\int^{\psi }\left[ \frac{%
f_{4}\left( \phi \right) }{2}+b\left( \phi \right) \right] d\phi }d\psi },
\label{c7a}
\end{equation}%
where $C_{12}$ is an arbitrary constant. By substituting Eq.~(\ref{c7a})
into the Riccati Eq.~(\ref{1}), the latter can be expressed as
\begin{equation}
\frac{dy}{dx}=a\left( x\right) +b(x)y+\frac{f_{4}\left( x\right) e^{-\int^{x}%
\left[ \frac{f_{4}\left( \phi \right) }{2}+b\left( \phi \right) \right]
d\phi }}{C_{12}+2\int^{x}a\left( \psi \right) e^{-\int^{\psi }\left[ \frac{%
f_{4}\left( \phi \right) }{2}+b\left( \phi \right) \right] d\phi }d\psi }%
y^{2}.  \label{TH1a}
\end{equation}

Therefore we obtain the following:

\textbf{Theorem 9}. If the coefficient $c(x)$ of the Riccati Eq.~(\ref{1})
satisfies the integral condition (\ref{c7a}), then the general solution of
the Riccati Eq.~(\ref{TH1a}) is given by

\begin{eqnarray}
y\left( x\right) &=&\frac{e^{\int^{x}\left[ f_{4}\left( \phi \right)
+b\left( \phi \right) \right] d\phi }}{C_{13}-\int^{x}\left\{ \frac{%
f_{4}\left( \psi \right) e^{-\int^{\psi }\left[ \frac{f_{4}\left( \phi
\right) }{2}+b\left( \phi \right) \right] d\phi }}{C_{12}+2\int^{\psi
}a\left( \omega \right) e^{-\int^{\omega }\left[ \frac{f_{4}\left( \phi
\right) }{2}+b\left( \phi \right) \right] d\phi }d\omega }\right\}
e^{\int^{\psi }\left[ f_{4}\left( \phi \right) +b\left( \phi \right) \right]
d\phi }d\psi }  \notag \\
&&+\frac{1}{2}\left\{ C_{12}+2\int^{x}a\left( \omega \right)
e^{-\int^{\omega }\left[ \frac{f_{4}\left( \phi \right) }{2}+b\left( \phi
\right) \right] d\phi }d\omega \right\} e^{\int^{x}\left[ \frac{f_{4}\left(
\phi \right) }{2}+b\left( \phi \right) \right] d\phi },
\end{eqnarray}%
where $C_{13}$ is an arbitrary constant of integration.

\subsection{Case 10: $\ a\left( x\right) =\frac{b^{2}\left( x\right)
-4c^{2}\left( x\right) f_{5}^{2}\left( x\right) }{4c\left( x\right) }+\frac{d%
}{dx}\left( -\frac{b}{2c}\pm f_{5}\right) $}

We assume that the coefficient $a\left( x\right) $ of the Riccati Eq.~(\ref%
{1}) satisfies the differential condition%
\begin{equation}
a\left( x\right) =\frac{b^{2}\left( x\right) -4c^{2}\left( x\right)
f_{5}^{2}\left( x\right) }{4c\left( x\right) }+\frac{d}{dx}\left( -\frac{b}{%
2c}\pm f_{5}\right) ,  \label{d9}
\end{equation}
where we have introduced an arbitrary function $f_{5}\left( x\right) \in
C^{\infty }(I)$ defined on a real interval $I\subseteq \Re $. By
substituting Eq.~(\ref{d9}) into the Riccati Eq.~(\ref{1}), the latter can
be expressed as

\begin{equation}
\frac{dy_{\pm }}{dx}=\left[ \frac{b^{2}\left( x\right) -4c^{2}\left(
x\right) f_{5}^{2}\left( x\right) }{4c\left( x\right) }+\frac{d}{dx}\left( -%
\frac{b}{2c}\pm f_{5}\right) \right] +b(x)y_{\pm }+c(x)y_{\pm }^{2}.
\label{TH13}
\end{equation}

Therefore we obtain the following:

\textbf{Theorem 10}. If the coefficient $a(x)$ of the Riccati Eq.~(\ref{1})
satisfies the differential condition~(\ref{d9}), then the general solutions
of the Riccati Eq.~(\ref{TH13}) are given by

\begin{equation}
y_{\pm }(x)=\frac{e^{\pm 2\int^{x}c\left( \phi \right) f_{5}\left( \phi
\right) d\phi }}{C_{\pm 14}-\int^{x}c\left( \psi \right) e^{\pm 2\int^{\psi
}c\left( \phi \right) f_{5}\left( \phi \right) d\phi }d\psi }-\frac{b\left(
x\right) }{2c\left( x\right) }\pm f_{5}\left( x\right) ,
\end{equation}

where $C_{\pm 14}$ are arbitrary constants of integration.

\section{Applications in physics: anisotropic stars in general relativity}

\label{sect3}

In standard coordinates $x^{i}=\left( t,r,\chi ,\phi \right) $ the line
element for a static spherically symmetric space-time takes the form \cite%
{HM}
\begin{equation}
ds^{2}=A^{2}\left( r\right) dt^{2}-V^{-1}\left( r\right) dr^{2}-r^{2}\left(
d\chi ^{2}+\sin ^{2}\chi d\phi ^{2}\right) .  \label{r1}
\end{equation}

For the metric (\ref{r1}), Einstein's gravitational field equations (where
natural units $8\pi G=c=1$ have been used throughout) describing the
evolution of the star take the form \cite{HM}
\begin{eqnarray}
\rho \left( x\right) &=&\frac{1-V\left( x\right) }{x}-2\frac{dV}{dx}, \\
p_{r}\left( x\right) &=&4V\left( x\right) \frac{d\ln A}{dx}+\frac{V\left(
x\right) -1}{x},  \label{r2}
\end{eqnarray}%
and
\begin{equation}
\left[ 1-2x\eta \left( x\right) \right] \frac{d^{2}A}{dx^{2}}-\left[ x\frac{%
d\eta }{dx}+\eta \left( x\right) \right] \frac{dA}{dx}-\left[ \frac{1}{2}%
\frac{d\eta }{dx}+\frac{\Delta \left( x\right) }{4x}\right] A=0,  \label{r3}
\end{equation}%
where we have used the coordinate transformation $x=r^{2}$, $\rho \left(
x\right) $ is the energy density, $p_{r}\left( x\right) $ is the pressure in
the direction $\chi ^{i}$ (normal pressure), and $p_{\perp }\left( x\right) $
is the pressure orthogonal to $\chi _{i}$ (transversal pressure),
respectively. $\chi ^{i}$ is the unit space-like vector defined as $\chi
^{i}=\sqrt{V}\delta _{1}^{i}$. We assume that $p_{r}\left( x\right) \neq
p_{\perp }\left( x\right) $.

The anisotropy parameter $\Delta \left( x\right) $ is defined as $\Delta
\left( x\right) =p_{\perp }\left( x\right) -p_{r}\left( x\right) $. We have
also defined the following functions,
\begin{equation}
V\left( x\right) =1-2x\eta \left( x\right) , \qquad \eta \left( r\right) =%
\frac{m\left( r\right) }{r^{3}},
\end{equation}%
and
\begin{equation}
2m\left( r\right) =\int_{0}^{r}\xi ^{2}\rho \left( \xi \right) d\xi ,
\end{equation}%
respectively. The function $m\left( r\right) $ represents the mass
distribution within radius $r$. By introducing a new function $u\left(
x\right) $, defined as%
\begin{equation}
u\left( x\right) =\frac{1}{A}\frac{dA}{dx},  \label{q1}
\end{equation}%
provides the metric function $A\left( x\right) $ in the form
\begin{equation}
A\left( x\right) =A_{0}e^{\int^{x}u\left( \phi \right) d\phi },  \label{q2}
\end{equation}%
with $A_{0}$ an arbitrary constant of integration. By substituting Eq.~(\ref%
{q1}) into Eq.~(\ref{r3}), the latter becomes a Riccati equation of the form
\begin{equation}
\frac{du}{dx}=\frac{\frac{1}{2}\frac{d\eta }{dx}+\frac{\Delta \left(
x\right) }{4x}}{1-2x\eta \left( x\right) }+\frac{x\frac{d\eta }{dx}+\eta
\left( x\right) }{1-2x\eta \left( x\right) }u-u^{2}.  \label{r4}
\end{equation}

The physical quantities $\eta \left( x\right) $ and $\Delta \left( x\right) $
are free arbitrary functions here. Now, by comparing the Riccati Eq.~(\ref%
{r4}) with the Riccati Eq. (\ref{1}), one gets the following settings
\begin{equation}
y\rightarrow u,\qquad a\left( x\right) \rightarrow \frac{\frac{1}{2}\frac{%
d\eta }{dx}+\frac{\Delta \left( x\right) }{4x}}{1-2x\eta },\qquad b\left(
x\right) \rightarrow \frac{x\frac{d\eta }{dx}+\eta }{1-2x\eta },\qquad
c\left( x\right) \rightarrow -1,
\end{equation}
respectively. By inserting these settings into the ten theorems presented in
Section II, we can obtain ten solutions of the Riccati Eq.~(\ref{r4}), and,
consequently, ten solutions of the interior Einstein's gravitational field
equations. In order not to have repetitive calculations, we shall not
present these results in the present paper.

\section{Conclusions}

\label{sect4}

In the present paper, we have obtained a number of integrability conditions
of the Riccati Eq.~(\ref{1}). All the theorems presented in this paper have
been obtained with the help of Eq.~(\ref{GS}). By using Eq.~(\ref{GP}), Eq.~(%
\ref{GS}) can be expressed equivalently as
\begin{equation}
y\left[ a\left( x\right) ,c\left( x\right) ,y_{p}\left( x\right) \right] =%
\frac{y_{p}\left( x\right) e^{\int^{x}\left[ c\left( \phi \right)
y_{p}\left( \phi \right) -\frac{a\left( \phi \right) }{y_{p}\left( \phi
\right) }\right] d\phi }}{C-\int^{x}c\left( \psi \right) y_{p}\left( \psi
\right) e^{\int^{\psi }\left[ c\left( \phi \right) y_{p}\left( \phi \right) -%
\frac{a\left( \phi \right) }{y_{p}\left( \phi \right) }\right] d\phi }d\psi }%
+y_{p}\left( x\right) ,  \label{GS1}
\end{equation}
and
\begin{equation}
y\left[ a\left( x\right) ,b\left( x\right) ,y_{p}\left( x\right) \right] =%
\frac{y_{p}^{2}\left( x\right) e^{-\int^{x}\left[ b\left( \phi \right) +%
\frac{2a\left( \phi \right) }{y_{p}\left( \phi \right) }\right] d\phi }}{%
C-\int^{x}\left[ \frac{dy_{p}\left( \psi \right) }{d\psi }-a\left( \psi
\right) -b\left( \psi \right) y_{p}\left( \psi \right) \right]
e^{-\int^{\psi }\left[ b\left( \phi \right) +\frac{2a\left( \phi \right) }{%
y_{p}\left( \phi \right) }\right] d\phi }d\psi }+y_{p}\left( x\right) .
\label{GS2}
\end{equation}
Note that the complexity of the general solutions Eqs.~(\ref{GS}), (\ref{GS1}%
), and (\ref{GS2}) depends on the coefficients $b\left( x\right) $ and $%
c\left( x\right) $, $a\left( x\right) $ and $c\left( x\right) $, $a\left(
x\right) $ and $b\left( x\right) $ respectively, and on the particular
solution $y_{p}\left( x\right) $, satisfying the Riccati Eq.~(\ref{1}).
Therefore, one may choose some simple particular solutions to obtain the
theorems in the present paper.

We have also discussed a physical application of the presented results. The study of general relativistic compact objects is of
fundamental importance for astrophysics. After the discovery of the pulsars and the explanation of their properties by assuming that they are rapidly rotating neutron stars, the theoretical investigation of superdense stars  has been done using both numerical and analytical methods and the parameters of neutron stars have been worked out by using a general relativistic approach.

Since the Einstein's gravitational field equations describing the interior
of compact stellar objects can be reduced to a Riccati equation, the
mathematical results obtained in the present paper can be used to model
massive stars, such as neutron stars and pulsars. The results may be useful whenever there is some
uncertainty regarding the actual equation of state of the
dense matter inside the general relativistic star. Note, however, that in order
to have physical stellar models, the interior solution for static fluid
spheres of Einstein's gravitational field equations must satisfy specific
general physical requirements. The following conditions have been generally
recognized to be crucial for anisotropic fluid spheres: (i) the density $%
\rho $, the radial pressure $p_{r}$ and the tangential
pressure $p_{\perp }$ should be positive inside the star; (ii) the
gradients $d\rho /dr$, $dp_{r}/dr$ and $dp_{\perp }/dr$ should be negative; (iii) inside the static configuration the
speed of sound should be less than the speed of light, i.e. $0\leq
dp_{r}/d\rho \leq 1$, $0\leq dp_{\perp }/d\rho \leq 1$;
(iv) the interior metric should be joined continuously with the exterior
Schwarzschild metric and; (v) the radial pressure $p_{r}$ must
vanish at the boundary $r=R$ of the sphere, but the tangential
pressure $p_{\perp }$ may not vanish. The physical interpretation
of the structure and properties of the compact relativistic star models
obtained by using the previous theorems as applied to the Riccati Eq.~(\ref%
{r4}) will be presented in a future publication.

\section{Acknowledgement}

We would like to thank to the anonymous referee for comments and suggestions that helped us to improve our manuscript. FSNL acknowledges financial support of the Funda\c{c}\~{a}o para a Ci\^{e}ncia e Tecnologia through the grants  CERN/FP/123615/2011 and CERN/FP/123618/2011.

\end{document}